\documentclass[11pt]{amsart}
\usepackage{amssymb}
\newcommand{\rfp}{{\mathcal{R}_+}}

\newcommand{\FT}{{\mathcal{F}}}
\newcommand{\h}{{\frac 12}}
\newcommand{\tu}{{\widetilde{V_0}}}
\newcommand{\wtp}{{\widetilde{P}}}
\newcommand{\mcl}{\mathcal{L}}
\def \intx {\stackrel{\,\circ}{X}}

\newcommand{\mcc}{\mathcal{C}}

\newcommand{\pa}{{\partial}}
\newcommand{\eps}{{\epsilon}}

\newcommand{\RR}{{\mathbb{R}}}

\newcommand{\hyp}{{\mathbb{H}}}

\newcommand{\lag}{{\mathcal{L}}}
\newcommand{\restrictedto}{{\upharpoonright}}
\newcommand{\WF}{{\operatorname{WF}}}

\newcommand{\OG}{{\mathsf{OG}}}
\newcommand{\CI}{{\mathcal{C}^\infty}}
\newcommand{\abs}[1]{{\left\lvert{#1}\right\rvert}}
\newcommand{\norm}[1]{{\left\lVert{#1}\right\rVert}}

\newcommand{\Lap}{{\Delta}}

\DeclareMathOperator{\sgn}{{sgn}}
\DeclareMathOperator{\sech}{{sech}}

\theoremstyle{plain}
\newtheorem{theorem}{Theorem}[section]
\newtheorem{proposition}[theorem]{Proposition}

\theoremstyle{definition}

\theoremstyle{remark}

\newtheorem{example}{Example}

\numberwithin{equation}{section}

\def\antonio/{{S\'{a} Barreto}}

\thanks{The first author acknowledges support from NSF grant DMS-0140657 and the
second author from grant DMS-0323021.}

\title{The radiation field is a Fourier integral operator}
\author{Ant\^{o}nio S\'{a} Barreto}
\address{Department of Mathematics\\
Purdue University\\
150 North University Street\\
West Lafayette IN 47907}
\author{Jared Wunsch}
\address{Department of Mathematics\\
Northwestern University\\
2033 Sheridan Rd.\\
Evanston IL 60208}

\begin{document}
\input epsf
\maketitle

\section{Introduction}

In this note, we exhibit explicitly the form of the ``radiation field'' of
F.~G.~Friedlander on two different types of manifolds: scattering
manifolds, and asymptotically hyperbolic manifolds.  The former class
consists of manifolds with ends that look asymptotically like the large
ends of cones, and includes a large class of asymptotically Euclidean
spaces, while the latter consists of spaces that resemble the hyperbolic
space at infinity, and includes quotients of hyperbolic space by certain
groups of motion.  In both cases we assume that there are no trapped
geodesics.  The radiation field is a measurement of the (rescaled)
asymptotic behavior of solutions to the wave equation, viewed from the
point of view of a rescaled time coordinate which in the asymptotically
Euclidean setting is simply $s=t-r,$ and restricted to the sphere at
infinity.  In particular, in $\RR^n,$ we define
$$
\rfp (s,\theta,z') = \lim_{r\to\infty} r^{(n-1)/2} V(s+r,r\theta,z')
$$ where $V(t,z,z')=(\cos t\sqrt \Lap, \sin t\sqrt\Lap/\sqrt\Lap)$ is the
solution operator to the wave equation.  Friedlander showed in
\cite{Friedlander3} that $\rfp$ is in fact a a translation representation
of the wave group in the sense introduced by Lax and Phillips
\cite{Lax-Phillips1}.  (S\'{a} Barreto \cite{SaBarreto1,SaBarreto2}
subsequently showed the unitarity of this map.)

In this paper, we show that the radiation field $\rfp$ on a manifold $X$
which is either a scattering manifold or an asymptotically hyperbolic
manifold, has as its Schwartz kernel a Lagrangian distribution, associated
to the conic Lagrangian defined by the graph of a ``sojourn relation''
relating points in $T^* \intx$ to points in $T^*(\RR\times \pa X),$ where
$\pa X$ is the boundary at infinity.  In the simple example where $\intx$
is just Euclidean space, the graph of the sojourn relation maps $T^* \RR^n$
to $T^* (\RR \times S^{n-1})$ roughly as follows: given $(z,\hat\zeta) \in
S^*(\RR^n)$, let $\gamma$ be the unique unit speed geodesic passing through
it.  We map $(z,\zeta)$ to the base point $(s,\theta)\in\RR\times S^{n-1}$
where $s$ is given by the ``sojourn time'' or the limit of $t-r$ along the
geodesic, and $\theta$ is the asymptotic direction in $S^{n-1}.$ The fiber
variables then measure the angle of contact the geodesic makes with
$S^{n-1}$ in a way made precise below.  This sojourn relation on scattering
manifolds was previously investigated by Hassell-Wunsch
\cite{Hassell-Wunsch1} in the context of the fundamental solution to the
Schr\"odinger equation on scattering manifolds; it is also closely related
to the sojourn time defined by Guillemin in \cite{Gu} in the study of the
high frequency asymptotics of the scattering matrix.  Note that the sojourn
time $s$ is more or less just the ``Busemann function'' used in
differential geometry.

In the special case in which we locally have a finitely many geodesics
$\gamma_n(t)$ beginning at a point $z\in \intx$ with asymptotic direction
$\theta$, the construction is simpler: $\rfp$ is a conormal distribution
with respect to the hypersurfaces $s=S_n(z,\theta)\equiv\lim
t-r(\gamma_n(t)).$ In this case, Fourier transforming yields a simple
result about the high-frequency asymptotics of the scattering Poisson
operator (better known, in the asymptotically hyperbolic case, as the
Eisenstein function).  This result is a weak generalization of a result of
Guillemin \cite{Gu}, who proved on compactly supported perturbations of
Euclidean space that not only do the scattering operator and Poisson
operator have the form discussed here, but a composition of FIO's gives the
asymptotics of the scattering matrix as well.  We are unable to perform
this composition owing to the local nature of our results in $z$, the
variable in $\RR^n$ (the location of the initial pole of the fundamental
solution).  Our results are also weaker in the sense that we obtain only
\emph{distributional} asymptotics of the Poisson operator, i.e.\ we must
mollify by convolution with the inverse Fourier transform of a compactly
supported cutoff in order to describe the asymptotics.  On the other hand
the results presented here are novel insofar as we do permit global
perturbations of the metric and folded sojourn relations, and we treat the
asymptotically hyperbolic case as well.

In the case of obstacle, rather than metric, scattering, results relating
the scattering matrix and the sojourn time were initially obtained by Majda
\cite{Majda1}, and used in the solution of inverse problems.  For further
applications of sojourn-time methods in inverse obstacle scattering, see
the survey by Petkov and Stoyanov \cite{Petkov-Stoyanov2}.  Similar results
to ours in the case of semiclassical scattering on $\RR^n$ have also been
obtained by Robert and Tamura \cite{Robert-Tamura1}; more recently,
Alexandrova \cite{Alexandrova} has extended these results to show that even
if the Hamilton flow is degenerate, the scattering matrix is a
semiclassical FIO.

We now discuss simple examples in which the Poisson operator is explicitly
known and the appearance of the sojourn time in the high-frequency behavior
is clear.

\begin{example}
On $\RR^n,$ the kernel of the scattering-theoretic Poisson operator,
evaluated at $(z,\theta)\in \RR\times S^{n-1},$ is just
$$
\left(\frac{i\lambda}{2\pi}\right)^{(n-1)/2} e^{i\lambda \theta \cdot z}
$$ (see \cite{Melrose44}); this is the operator mapping high-frequency
incoming scattering data to a generalized eigenfunction with eigenvalue
$\lambda^2$ (see \S\ref{subsec:poisson} for a precise definition.).

On the other hand, there is a unique geodesic beginning at $z$ with
asymptotic direction $\theta:$ it is just $z+t\theta.$ The sojourn time
along the geodesic is
$$
\lim\abs{t-z+t\theta}=-\theta \cdot z,
$$ hence exactly the phase of the adjoint of the Poisson operator.
\end{example}

\begin{example}
The kernel of the Poisson operator on hyperbolic space is better known as
the ``Eisenstein function;'' it is given by
$$
E(\frac n2 +i\lambda,y,z)=\lim_{x\to 0^+} x^{-n/2-i\lambda}R\left(\frac n2+i\lambda,x,y,z\right)
$$ where $R(\lambda)$ is the resolvent, normalized to be
$(\Lap-\lambda^2-\frac{n^2}4)^{-1}$ and where we work in the usual
coordinates on the half-space with defining function $x.$  On $\hyp^3$,
the resolvent is just (see equation (6.8) of \cite{Mazzeo-Melrose1}):
$$
R(\frac n2+i\lambda) = C \frac{e^{-i\lambda \delta}}{\sinh \delta},
$$
where $\delta$ is the hyperbolic distance.  The phase is thus $-\delta,$
which is asymptotic to
$$
-\log \frac{x^2+ (x')^2+ (y-y')^2}{xx'},
$$
hence, switching primed and unprimed variables, the phase of
$E(\frac n2+i\lambda,y',z)$ is just
$$
\phi(x,y,y')\equiv-\log \frac{x^2+(y-y')^2}{x}.
$$

Now we compare this phase to the sojourn time.  Given a point $z\in \hyp^3$
and $y' \in \pa \hyp^3,$ there is a unique geodesic $\gamma(t)$ starting at
$z,$ and approaching $y'$ as $t\to \infty$ (this is the analogue of the
``asymptotic direction'' $\theta \in S^{n-1}$ in the Euclidean case). We
now define the sojourn time for this geodesic as
$$
S=S(z,y')=\lim t+\log x(\gamma(t))
$$
Using translation invariance in the boundary variables, it suffices to
compute with $y'=0.$ A unit speed geodesic from an arbitrary point
$(x,y)\in \hyp^3$ to $(0,0) \in \pa \hyp^3$ can be parametrized as
$$
x=C\sech(t+t_0),\ y= C \tanh(t+t_0)+D,
$$
whence we compute $\phi(x,y,0)=-\log 2-\log C+t_0.$  On the other hand, the
sojourn time along such a geodesic is
$$
S(x,y,0)\equiv\lim_{t\to +\infty} t+ \log x(t)=\log 2 +\log C-t_0,
$$
hence agrees with minus the phase $\phi(x,y,0).$
\end{example}

\section{Scattering manifolds}

\subsection{Radiation field}\label{subsec:radfield}
Let $X$ be a $\mcc^\infty$ compact manifold with a
boundary.   A scattering
metric (defined originally by Melrose \cite{MR95k:58168}) is a metric on a
manifold with boundary that can be brought to the form
\begin{equation}
g=\frac{dx^2}{x^4}+ \frac h{x^2}
\label{scmetric}
\end{equation}
with $x$ a boundary defining function and $h$ a smooth tensor that
restricts to $x=0$ to give a metric $h_0$ on $\pa X.$ This form is modeled after
the metric on asymptotically Euclidean space, radially compactified.  It
was shown by Joshi-S\'{a} Barreto \cite{Joshi-SaBarreto2} that there is a
normal form for a scattering metric: in which a neighborhood of the
boundary admits a product decomposition $[0,\epsilon) \times \pa X$ with
local coordinates $(x,y)$ in which $g$ takes the form \eqref{scmetric} with
$h=h(x,y,dy)$ a smooth family in $x$ of metrics on $\pa X.$ Note that in
these product coordinates, the rays $y=\text{constant}$ are (infinitely
extended) geodesics.  We henceforth assume that our metric is in normal
form.  We further make the geometric assumption that there are no trapped
rays in $\intx.$

Following Friedlander \cite{Friedlander3} we define the \emph{forward
  radiation field} on a scattering manifolds by
$$
\rfp(f_1,f_2)(s,y) = (x^{-(n-1)/2} D_t Hu)(s+x^{-1},x,y)\restrictedto_{x=0}
$$ where $H=H(t)$ is the Heaviside function and $u=u(t,x,y)$ is the solution to
\begin{gather}
\Box u=0, \quad (u,D_t u) \restrictedto_{t=0} = (f_1,f_2).\label{we1}
\end{gather}
Friedlander  showed that $\rfp(f_1,f_2)\in
\CI(\RR\times \pa X)$ provided $f_1,f_2$ are smooth and compactly supported
in $\intx.$  Thus the Schwartz kernel
$$
\rfp(s,y,z)
$$ is defined on $\RR\times \pa X;$ here we have used coordinates $(s,y,z)$
on $\RR\times \pa X \times \intx.$ First we will show that the radiation
field is an FIO and then, under some nondegeneracy conditions, we will compute
 its symbol.  We use coordinates on $T^*(\RR\times \pa X \times \intx)$
 defined by the canonical one-form
$$
s \, d\sigma + y \, d\eta + z \, d\zeta,
$$
and will employ the notation $\hat \zeta=\zeta/\abs\zeta.$
\begin{theorem}\label{propeu}
Let $X$ be a nontrapping scattering manifold.  Then
$$
\rfp \in (  I^1 (\RR\times \pa X\times \intx, \Lambda';
\Omega^{1/2}),\ I^0 (\RR\times \pa
  X\times \intx, \Lambda'; \Omega^{1/2})), 
$$ where $\Omega$ is the density bundle on $\RR\times \pa X \times \intx,$
and $\Lambda_\pm'$ is the conic Lagrangian manifold associated to the graph of the
``sojourn relation'':
\begin{multline*}
\Lambda_\pm=\Lambda_+ \cup \Lambda_-,\\
\Lambda_\pm=\bigg\{z,\zeta, s=\lim_{t\to \pm\infty} t-x^{-1}(\exp_z(t\hat
\zeta)), \sigma=\pm \abs{\zeta}, \\ y=\lim_{t\to\pm\infty} y(\exp_z(t\hat \zeta)),
\eta_i=\lim_{t\to\pm \infty} \pm \abs{\zeta}(h_0)_{ij}(dy^j/dx)(\exp_z(t\hat \zeta))\bigg\}.
\end{multline*}
\end{theorem}
\begin{proof} We know from the non-trapping assumption and from
Theorem 1.1 of \cite{Duistermaat-Guillemin1} that the kernel of $\cos
t\sqrt{\Delta},$ which is the solution to \eqref{we1} with $f_1=\delta(z)$
and $f_2=0,$ a sum of forward and backward parts (corresponding to a choice
of $\pm$ in our notation below), each of which is a Lagrangian distribution
in $I^{-\frac{1}{4}}\left( \RR \times X\times X; \mcc_\pm \right),$ where
\begin{gather}
\begin{gathered}
\mcc_{\pm}=\{  (t,z,z',\tau,\zeta,\zeta'): (z,\zeta), \;
 (z',\zeta') \in T^* X\setminus 0,
(t,\tau) \in T^* \RR \setminus 0, \\ \; \tau =\pm \sqrt{L(z,\zeta)}, \;
(z,\zeta)= \Phi^{\pm t}(z',\zeta') \}.
\end{gathered}\label{DG1}
\end{gather}
Here $L(z,\zeta)$ denotes the symbol of $\Delta$ and $\Phi^t(z',\zeta')$
denotes the flow along $H_L$ in $T^* X \setminus 0.$ We want to understand
what happens to the kernel of $\cos t\sqrt{\Delta}\;\ $ if one first makes
the change $s=t-1/x$ and then takes the limit $x\rightarrow 0.$ So
we need to understand the effect of these operations on $\mcc_{\pm}$ and
the distribution associated with it.

As in \cite{SaBarreto1}, let $P=x^{-2-(n-1)/2} \Box x^{(n-1)/2}$
and change variables, replacing $t$ by $s=t-1/x.$  We find as in 
\cite{SaBarreto1} that
$$
P=2 \frac{\pa}{\pa x}\frac{\pa}{\pa s}+ x^2\frac{\pa^2}{\pa x^2}-
\Lap_h+ A \frac{\pa}{\pa s} + (2x+x^2A)\frac{\pa}{\pa x} +
\left(\frac{n-1}{2}\right)\left(\frac{3-n}{2}+ xA\right) 
$$
with $A(x,y)=\pa_x \log\abs{h}^{1/2}$ and $\Lap_h$ is the nonnegative 
Laplacian on 
$\pa X$ with respect to the metric $h.$
The symbol of $P$ in these coordinates is given by
$$ p=-2\xi\sigma-x^2 \xi^2-h(x,y,\eta)
$$
and the Hamilton vector field by
$$
H_p = -2(\sigma+x^2\xi) \frac{\pa}{\pa x} - 2 \xi \frac{\pa}{\pa s} +
\left( 2 x \xi^2 + \frac{\pa h}{\pa x}\right) \frac{\pa}{\pa \xi} - H_h
$$ ($(z, \zeta)$, the coordinates in the right factor, are left
invariant by the flow).

 If $V_0=x^{-\frac{n-1}{2}}U_0,$ the equation
\eqref{we1} with initial data $f_1=\delta(z)$ and $f_2=0$ becomes
\begin{equation}
P V_0=0, \quad V_0\restrictedto_{s=-\frac{1}{x}} =
x^{-\frac{n-1}{2}}\delta(z),\quad 
D_t V_0\restrictedto_{s=-\frac{1}{x}}=0, \;\ x>0
\label{IVP1}
\end{equation}
 We remark 
that we are ignoring half-density factors, as they are irrelevant to this 
construction.  Notice that the operator $P$ extends to $x\leq 0$ as a 
strictly hyperbolic differential operator $\wtp$ and that
$s=-\frac{1}{x}$ is a space-like surface for $\wtp.$  One can think of this as 
being an extension of $P$ to the double manifold $X^2=(X\sqcup X)/ \pa X.$
 Therefore it follows from 
the existence of a fundamental solution to the Cauchy problem for
strictly hyperbolic operators, see for example Theorem 5.1.6 of \cite{Dui},
that the solution to \eqref{IVP1}, with $P$ replaced by $\wtp,$
is a Lagrangian distribution $\tu$ of class
\begin{gather}
\tu \in I^{-\frac{1}{4}}\left( \RR \times X\times \intx; \mcl\right), 
\label{lagreg}
\end{gather}
 where $\mathcal{L}$ denotes the Lagrangian in $T^*(\RR\times X \times
 \intx)$ obtained by flowing $ N^*\{z=z', \;\ t=0\} \cap \Sigma_p\subset
 T^* (\RR\times X\times \intx)$ along the integral curves of
 $H_{\widetilde{p}},$ where $\widetilde{p}$ is the principal symbol of
 $\wtp,$ and $\Sigma_{\widetilde{p}}$ denotes the characteristic variety of
 $\wtp.$ By the uniqueness of solutions to the Cauchy problem, the
 restriction of $\tu$ solution to $x>0$ is equal to $V_0,$ the solution to
 \eqref{IVP1}.

We remark that the extension $\wtp,$ and consequently the definition of
$\mcl$ in $\{x<0\},$ are not unique.  However the extensions to $\{x=0\}$
are.

Since $N^*(\{x=0\})\cap \mcl=0,$ $U\equiv \tu \restrictedto_{x=0}$ is a
Lagrangian distribution of order $1$ on $\RR\times \pa X \times \intx$ with
respect to $\lag\restrictedto_{x=0}.$ Hence for any compactly supported
distribution $u$, $\rfp ( u,0) = D_t \widetilde{U} u$ where $\widetilde{U}$
is a Lagrangian of order $0$ with respect to $\lag\restrictedto_{x=0}.$

To identify $\lag\restrictedto_{x=0}$ geometrically we observe that under
the flowout of $H_p$, $\sigma$ is conserved, hence by homogeneity we need
only consider $\sigma=\pm 1$ (note that $\sigma \neq 0$ on $\lag \setminus
0$).  Let $\lag_\pm$ denote the two components corresponding to different
signs of $\sigma.$ By definition, we certainly have $s=
\lim_{t\to\pm\infty} (t-x^{-1}),$ and $y$ is the limiting location of
geodesic flow (forward or backward according to $\sgn \sigma$) in $\pa X.$
Furthermore, setting $\sigma=\pm 1$ we have $dy^i/dx =
(-2h^{ij}\eta_j)/(-2(\pm 1+x^2 \xi)),$ which approaches $\pm h^{ij}\eta_j$
as $x\to 0.$ Thus the restriction of $\lag_\pm$ equals $\Lambda_\pm$ as
defined in the statement of the theorem.

Similarly, solving
\begin{equation}
P U_1=0, \quad U_1\restrictedto_{t=0} =0,\quad 
D_t U_0\restrictedto_{t=0}=\delta(z)
\label{IVP2}
\end{equation}
gives a Lagrangian of order $-1$ when restricted to the boundary.  Hence
for any distributions $f_0,f_1,$
$$
\rfp(f_0,f_2)=(D_t U_0 f_0, D_t U_1 f_1) \restrictedto_{x=0}
$$
is a FIO of the asserted kind.
\end{proof}

It remains to calculate the symbol of $\rfp$ and for this we need some
extra assumptions.  So we suppose further that for all $z$ contained in an
open set $U_1$ in $\intx,$ all $y$ contained in an open set $U_2$ in $\pa
X$, there exist a finite number of unit speed geodesics
$\gamma_{n}(z,y,t),$ $n=1,\dots, N$ such that
$$
\gamma_{n}(z,y,0)=z,\quad \lim_{t\to +\infty} \gamma_{n}(z,y,t) = y,
$$
and such that the transformation
$$
\frac{\pa^2 \gamma_{n}}{\pa y\pa t} \restrictedto_{t=0}
$$
mapping
$$
T(\pa X) \to T(S(\intx))
$$ is invertible for all $n,$ and $y\in U_2,$ $z\in U_1.$
Subject to these assumptions, the Lagrangian $\Lambda$ is projectable
onto the $z,y$ variables.  Letting $$S_{n}(z,y)= \lim_{t\to
  +\infty} t-x^{-1}(\gamma_{n}(z,y,t))$$ (the ``sojourn times''), we
then find that $\rfp(0, \delta_z)$ is conormal to the surfaces
$s=S_{n}(z,y)$ in $\RR_s \times U_1 \times U_2.$
\begin{theorem}\label{propeusy}
Subject to the nondegeneracy assumptions above, the symbol of $\rfp,$
evaluated at $(z,y,\sigma) \in N^* \{s=S_{n}\}$ equals
\begin{equation}\label{rsymbol}
\h i^{k_n}\abs{\frac{\pa y^{(n)}}{\pa \hat\zeta}}^{-1/2}\abs{ds\,
  dh_y \, dg_{z}}^{1/2}\cdot (\sigma^{(n+1)/2},
  \sigma^{(n-1)/2}).
\end{equation}
where $k_n$ is the number of conjugate points encountered by
$\gamma_{n}(z,y,t)$ with $t \in (0,\infty)$ and $\abs{{\pa
y^{(n)}(z,y)}/{\pa \hat\zeta}}$ is short for the Jacobian of the map
$$
S^*_z(\intx) \to \pa X
$$ given by the limit of geodesic flow, evaluated at the initial
codirection of $\gamma_n.$
\end{theorem}
Note that we have written the Jacobian factor in the above form to
emphasize the analogy with the differential scattering cross-section in
Guillemin's results \cite{Gu}.
\begin{proof} 
Let $\OG\subset S^* (\intx)$ denote the outgoing set, i.e.\ the set on
which $dx/dt<0$ along the bicharacteristic flow.  Owing to our nontrapping
assumption, the cosphere bundle of a compact set $S^* (K)\subset S^*
(\intx)$ eventually maps into $\OG$ under the bicharacteristic flow after
time $T\gg 0.$ Let us fix such $K$ and $T,$ with $K$ chosen such that our
nondegeneracy assumption holds for all geodesics beginning in $K$ with
limits some open set in $y.$

Net $Z$ be a compactly supported pseudodifferential operator with
$\WF'Z\subset \OG,$ $WF'(1-Z) \cap \OG'=\emptyset,$ with $\OG'$ a conic
subset of $\OG$ with compact projection, chosen so that the flowout of $S^*
K$ for time $T$ lies inside $\OG'.$ 

Let $Z_j'$ be a microlocal partition of unity over $K$ such that for all
$j,$ and all $z\in \intx,$ $\WF' Z_j'\cap \pi^{-1}(z)$ contains at most a
single point such that the geodesic emanating from this point ends up at $y
\in \pa X.$ By our hypotheses on the nondegeneracy of geodesics, we may
further choose $Z$ supported sufficiently close to $y$ that if
$(z',\zeta')$ and $(z,\zeta)$ are canonical coordinates on the cotangent
bundles of left and right factors, then $(\zeta', z)$ are coordinates on
$$
\Lambda_{\pm} \cap \pi_L^* \WF' Z \cap \pi_R^* \WF' Z_j'.
$$ We may further arrange, by working sufficiently close to $\pa X,$ that
there are no conjugate points for bicharacteristics beginning on $\WF' Z:$
examination of the Hamilton flow in ``scattering coordinates'' shows that
the tangent vectors to geodesics emanating from $\WF' Z$ approach $-x^2 \pa
x,$ and while the sectional curvature of a scattering manifold is $O(x^2)$,
the sectional curvature of a plane containing $x^2 \pa x$ is $O(x^3).$
Hence certainly the sectional curvature along a plane containing the
tangent to the geodesic is $o(x^2).$ By a simple variant on the Rauch
comparison theorem, this is sufficient to ensure the absence of conjugate
points.  (See \cite{Melrose43} for a description of the bicharacteristic
flow in scattering coordinates and, for instance Theorem 4.5.1 of
\cite{Jost} for an account of the relevant comparison theorem.)

Now let $W(t)$ denote the propagator for the wave group.
The symbol of $W(T) Z_j'$ is $$\h\sigma(Z_j')(z,\zeta) i^{k_n} \abs{dz\,
d\zeta}^{1/2}  \begin{pmatrix} 1 & \abs\zeta^{-1}_g
  \\ \abs{\zeta}_g & 1 \end{pmatrix}$$
where $k_n$ is the number of conjugate
 points encountered.
Changing to coordinates $(\zeta', z)$ gives
$$
\h \pi_R^* \sigma(Z_j') i^{k_n} \abs{\frac{\pa \zeta'}{\pa \zeta}}^{-1/2}
\abs{d\zeta' \, dz}^{1/2} \begin{pmatrix} 1 & \abs{\zeta'}^{-1}_g
  \\ \abs{\zeta'}_g & 1 \end{pmatrix}.
$$ By invariance under the bicharacteristic flow (see for instance
Proposition 4.3.1 of \cite{Dui}), the symbol of $\rfp Z$ is
$$
\h \pi_R^* \sigma(Z) 
\abs{\frac{\pa(s,y)}{\pa{\zeta'}}}^{-1/2} \abs{ds\, dy\,
  dz'}^{1/2}(\abs{\sigma}, 1).
$$ We split the coordinate $\zeta'$ into $\abs\zeta', \hat\zeta'$ (where the
latter should really be regarded as $n-1$ components of $\hat \zeta'$).  By
homogeneity of the flow we have $\pa s/\pa \abs\zeta'=1,$ while $\pa y/\pa
\abs \zeta'=0,$ hence the symbol of $\rfp Z$ can in fact be written
$$
\h \pi_R^* \sigma(Z) 
\abs{\frac{\pa(y)}{\pa{\hat\zeta'}}}^{-1/2} \abs{ds\, dy\,
  dz'}^{1/2}(\abs{\sigma}, 1).
$$

The difference
$$
\rfp(s) Z_j'  - \rfp(s-T) Z W(T) Z_j'
$$ is a smoothing operator.  Hence applying the calculus of FIO's to the
above results and patching together the partition of unity $Z_j',$ we find
that the symbol of $\rfp$ is given globally by \eqref{rsymbol}.
\end{proof}

\subsection{Poisson operator}
\label{subsec:poisson}
Let $P(\lambda)$ denote (the Schwartz kernel of) the \emph{Poisson
operator}, i.e.\ the operator such that for any $g \in \CI(\pa X)$ there
exists $u \in \CI(X)$ with $(\Lap-\lambda^2) u=0$ and
$$
u= e^{i\lambda/x} x^{(n-1)/2} g + e^{-i\lambda/x} x^{(n-1)/2} g_- + u'
$$
with $g' \in \CI (\pa X),$ $u'\in L^2(X;g).$

As a corollary of the results in the preceding section, we conclude the
following:
\begin{proposition}
Suppose as in the preceding proposition that there exist a finite number of
nondegenerate geodesics $\gamma_{n}$ from $z$ to $y$.  Let $\check \phi
\in \mathcal{C}_c^\infty (\RR).$ Then as $\abs{\lambda} \to \infty,$
$$ \phi(\lambda)*P(\lambda)^* (y,z) \sim\phi(\lambda)* \sum_{n=1}^N
i^{k_n}e^{i \lambda S_n} \left(\frac{\lambda}{2\pi i}\right)^{(n-1)/2}
\abs{\frac{\pa y^{(n)}}{\pa \hat\zeta}}^{-1/2}
$$
\label{prop:poisson}
\end{proposition}
\begin{proof}
We simply use the fact that
$$
P(\lambda)^*(y,z)=-2 \FT \rfp(0,\delta_z)
$$
(see \cite{SaBarreto1}).
\end{proof}

Note: if we knew more about energy decay, along the lines of having
good estimates for the energy norm
$$
E(t,\lambda) = \norm{V(t)f}_{H_E(t-\lambda)}=
\int_{t-\frac{1}{x}>\lambda} \left( |\nabla V(t) f|^2 + 
\left| \frac{\pa V(t) f}{\pa t}\right|^2 \right) \; d \operatorname{vol}_{g}
$$ with initial data $f$ compactly supported, we would be able to get
better estimates for the decay as $s\to +\infty$ of $\rfp,$ and hence drop
the mollifier $\phi$ from the statement of this proposition.  Friedlander
proves in \cite{MR82f:35112} that
$$
\lim_{\lambda \to \infty} \lim_{t\to\infty} E(t,\lambda)=0
$$ where $f$ is finite energy initial data; we would need a good deal
more however.

\section{Asymptotically hyperbolic manifolds}

An \emph{asymptotically hyperbolic}, or \emph{conformally compact}
manifold is a manifold $X$ with boundary equipped with a metric $g$ and
defining function $x\geq 0$ such that $x^2 g=H$ is a smooth metric on $X$,
nondegenerate at $\pa X.$  We further take $\abs{dx}_H=1$ on $\pa X,$
which ensures that the sectional curvatures approach $-1$ at $\pa X$ (see,
for instance, \cite{Mazzeo-Melrose1}.)

Notice that $H\restrictedto_{\pa X}$ is only defined by $g$ modulo a
conformal factor. It is shown in \cite{Joshi-SaBarreto1}, see also
\cite{graham}, that given a conformal representative of
$H\restrictedto_{\pa X}$ there exists a unique boundary defining function
$x$ such that
\begin{gather*}
g=\frac{dx^2 + h(x,y,dy)}{x^2} \;\ \text{ in } \;\ \pa X \times [0,\eps).
\end{gather*}
From now on fix these coordinates. The definition of the radiation fields will
depend on this choice of $x,$ or the conformal representative of 
$H\restrictedto_{\pa X}.$

We refer the reader to \cite{Joshi-SaBarreto1} for a discussion of the Eisenstein function
on an asymptotically hyperbolic manifold.

In this setting, we can prove more or less the same results as in the case
of scattering manifolds; the analogous results are as follows:

\begin{theorem}\label{thm:hyp}
Let $X$ be a nontrapping asymptotically hyperbolic manifold of dimension $n$.

\noindent (1) We have
$$
\rfp \in (  I^1 (\RR\times \pa X\times \intx, \Lambda';
\Omega^{1/2}),\ I^0 (\RR\times \pa
  X\times \intx, \Lambda'; \Omega^{1/2})), 
$$ where $\Lambda_\pm'$ is the conic Lagrangian manifold associated to the
graph of
\begin{multline}\label{ah:lag}
\Lambda_\pm=\Lambda_+ \cup \Lambda_-,\\ \Lambda_\pm=\bigg\{z,\zeta,
s=\lim_{t\to \pm\infty} t +\log x(\exp_z(t\hat \zeta)), \sigma=\pm
\abs{\zeta}, \\ y=\lim_{t\to\pm\infty} y(\exp_z(t\hat \zeta)),
\eta_i=\lim_{t\to\pm \infty} \pm x^{-1}
\abs{\zeta}(h_0)_{ij}(dy^j/dx)(\exp_z(t\hat \zeta))\bigg\}.
\end{multline}

\noindent (2)
Subject to the nondegeneracy assumptions of \S \ref{subsec:radfield}, the symbol of $\rfp$
is given by \eqref{rsymbol}.

\noindent (3)
Let $E(\frac n2 +i\lambda,y,z)$ denote the (transpose of the) Eisenstein function.  Let
$\check \phi \in \mathcal{C}_c^\infty (\RR).$ Subject to the nondegeneracy
assumptions of \S \ref{subsec:radfield}, as $\abs \lambda \to \infty,$
$$ \phi(\lambda)*E(\frac n2+ i \lambda,y,z) \sim\phi(\lambda)* \frac{i}{2\lambda}\sum_{n=1}^N
i^{k_n}e^{i \lambda S_n} \left(\frac{\lambda}{2\pi i}\right)^{(n-1)/2}
\abs{\frac{\pa y^{(n)}}{\pa \hat\zeta}}^{-1/2}
$$
\end{theorem}
\begin{proof}
It is shown in \antonio/ \cite{SaBarreto2} if $X$ has dimension $n,$ and
 $u=\cos\left( t \sqrt{\Delta}\right)$ 
than the rescaled fundamental solution, $v=x^{-(n-1)/2} u(s-\log x, x,y)$,
satisfies
\begin{gather}
\begin{gathered}
Pv \equiv \left( \frac \pa{\pa x}\left(2 \frac \pa {\pa s} + x \frac \pa {\pa
  x}\right) - x \Lap_h + A \frac \pa {\pa s}+ A x \frac \pa {\pa x} + \frac{n-1}{2}
 A \right) v=0  \\
v\restrictedto_{s=\log x}=x^{-(n-1)/2}\delta(z), \;\ 
\frac{\pa v}{\pa s}\restrictedto_{s=\log x}=0, \;\ x>0.
\end{gathered}\label{phyp}
\end{gather}
The symbol of $P$ 
is thus $-(2\xi \sigma + x \xi^2 + x h(x,y,\eta)),$ hence the flow
is
$$
H_p = -2(\sigma+ x\xi) \pa_x +(\xi^2+ h+ x\frac{\pa h}{\pa x})\pa_\xi- 
2\xi\pa_s- xH_h.
$$ As in the scattering case the operator $P$ is strictly hyperbolic in
$x>0$ and $\{s=\log x\}$ is space-like. Moreover $P$ has an extension
$\wtp$ to a neighborhood of $\{x=0\}.$ However, $\wtp$ is not strictly
hyperbolic at $x=0.$ So in principle we can only guarantee that $v$
satisfies \eqref{lagreg}, where $\mcl$ is defined with respect to $P$ in
\eqref{phyp}, when $x>0.$

Since $\wtp$ can be chosen to be of real principle type, the parametrix
construction of the Cauchy problem can be carried though across $\{x=0\}.$
This guarantees that in a neighborhood $W$ of $\pa X$ there exists
$$\tu \in I^{-\frac{1}{4}}\left( \RR \times W \times \intx; \mcl\right)$$
 such that
\begin{gather}
\begin{gathered}
P \tu =f \in \mcc^\infty( \RR \times W) \\
\tu\restrictedto_{s=\log x}-\delta(z)=g \in \mcc_0^\infty(\intx), \\
\frac{\pa \tu}{\pa s}\restrictedto_{s=\log x}=h\in \mcc_0^\infty(\intx).
\end{gathered}\label{cicp}
\end{gather}
Moreover by finite speed of propagation one has that $f$ is supported in
$s>s_0,$ for some $s_0.$
  Again using the hyperbolicity of $P$ in $x>0,$ 
there exists $V \in C^\infty\left(W\cap \{x>0\}\right)$  satisfying
\begin{gather*}
P V =f \;\ \text{ in } \;\ x>0 \\
V\restrictedto_{s=\log x}=g, \quad 
\frac{\pa V}{\pa s}\restrictedto_{s=\log x}=h,
\end{gather*}
with $V$ supported in $s>s_1.$ Since $f$ is supported in $s>s_0$ and is
smooth up to $\pa X,$ and the initial data is compactly supported, the
proof of Theorem 2.1 of \cite{SaBarreto2} (in particular, the extension of
the energy estimates to the inhomogeneous equation) shows that $V$ has a
smooth extension up to $\pa X.$ Therefore the solution to \eqref{phyp}
satisfies $v \in I^{-\frac{1}{4}}\left( \RR \times X \times \intx;
\mcl\right)$ up to $\pa X.$

The transversality of the flow to $\{x=0\}$ and the facts that $\sigma\neq
0$ on the characteristic variety for $x>0$ and $H_p\sigma=0$ imply that
$v\restrictedto_{x=0}$ is a Lagrangian distribution of order $0$ on
$\RR\times \pa X \times \intx$ with respect to the Lagrangian
$\lag\restrictedto_{x=0},$ where $\lag$ is the flowout of the lift of
$N^*\Delta \subset T^* (\intx \times \intx)$ to $\Sigma_p.$ We find in this
setting that $\lim_{x\to 0} x^{-1} dy^i/dx=\eta_j h_0^{ij}/\sigma,$ hence
$\lag\restrictedto_{x=0}$ has the form \eqref{ah:lag}.

The remainder of the proof is the same as in the scattering case, using the
 additional fact from \cite{SaBarreto2} that $E(n/2+i\lambda,y,z)=
 -(i/\lambda)\FT \rfp(0,\delta_z).$
\end{proof}

 \bibliography{all}
\bibliographystyle{amsplain}
\end{document}